\documentstyle[12pt]{article}
\begin{document}
%\baselineskip=24
\noindent $$\mbox{The Effect of Finite Memory Cutoff on Loop Erased Walk in} Z^{3}$$
\noindent $$\mbox{Wei-Shih Yang}$$
$$\mbox{Department of Mathematics}$$
$$\mbox{Temple University, Philadelphia, PA 19122}$$
$$\mbox{Aklilu Zeleke}$$
$$\mbox{Math and Computer Science Department}$$ 
$$\mbox{Alma College, Alma MI 48801}$$
\noindent {\bf Abstract} \noindent Let $S_{n}$ be a simple random walk (SRW) defined on $Z^{3}$. We
construct a stochastic process from $S_{n}$ by erasing loops of length at most
$N^{\alpha}$, where $\alpha \in (0, \infty)$ and $N$ is the scaling parameter that
will be taken to infinity in determining the limiting distribution. We call this
process the $N^{\alpha}$ loop erased walk ($N^{\alpha}$ LEW). Under some
assumptions we will prove that for $0 < \alpha  < \frac{1}{1+2\zeta}$, the
limiting distribution is Gaussian. Here $\zeta$ is the intersection
exponent of random walks in $Z^{3}.$ For $\alpha > 2$ the limiting distribution is equal to the
limiting distribution of the loop erased walk.\\
\vspace{3mm}


\vspace{300mm}

\noindent {\bf 1. INTRODUCTION}\\
\noindent A loop erased walk (LEW) is a stochastic process constructed from the simple random walk (SRW)
 by erasing paths that lead to the formation of loops. G. Lawler has proven that the limiting distribution 
of LEW in $Z^{d}, d \geq 4$ is Gaussian (see [L1]). The low dimensional cases
remain open. It is conjectured, however, that it is non-Gaussian.\\

\noindent In  this paper we consider a stochastic process constructed
from the SRW in $Z^{3}$ by erasing loops using only  finite memory. At
each step the first $N^{\alpha}$ loops will be erased (see  section 2 for
the definition of $N^{\alpha}$ loops). Here $\alpha \in  [0,\infty]$, and
$N$ is a scaling parameter which will be taken to $\infty$ in 
determining the limiting distribution of LEW. We call this
process the $N^{\alpha}$ LEW. Note that $\alpha=0$ is the
case of SRW and $\alpha=\infty$ is that of LEW. Under some 
assumptions we will prove that the $N^{\alpha}$ LEW has a
Gaussian distribution for $0 < \alpha < \frac{1}{1+2\zeta}$, where
$\zeta\;\mbox{is the intersection exponent of random walks
in}\;Z^{3}$. For $\alpha > 2$ we will show that the $N^{\alpha}$ LEW has
the same limiting  distribution as the original LEW. It can be implied
from our work that if there is a
%Our work lays down a foundation for the numerical
%study of the
%loop erased walk in three dimensions. It can be implied from our
critical point $\alpha_{c}$ then it must be between
$\frac{1}{1+2\zeta}\;\mbox{and}\;2$. The existence of $\alpha_{c}$ and
the behavior of  the $N^{\alpha}$ loop erased walk for $\frac{1}{1+2\zeta}
< \alpha < 2$  remain open.\\

\noindent{\bf 2. The $N^{\alpha}$ LOOP ERASED WALK}\\
\vspace{1mm}
\noindent Let $\lambda=[S_{i},S_{i+1},...S_{j}]$ be a segment
of a path of an SRW. We say that $\lambda$ forms an $N^{\alpha}$ loop if
$S_{i}=S_{j}\;\mbox{and}\;0\;\leq\;|i-j|\;\leq\;N^{\alpha}\;
\mbox{for some fixed}\;N\;\mbox{and}\;\alpha$. 
\noindent Let $\sigma_{\alpha}(0)=\mbox{sup}\{j: S(j)= 0, |j|\leq\;N^{\alpha}\}$, and
for $i>0\;\sigma_{\alpha}(i)=\mbox{sup}\{j>\sigma_{\alpha}(i-1):
S(j)=S(\sigma_{\alpha}(i-1)+1), |j-\sigma_{\alpha}(i-1)-1|\leq\;N^{\alpha}\}$. We define the 
$N^{\alpha}$ LEW by ${\hat{S}_{i}^{(N)}}=S(\sigma_{\alpha}(i))$. From now on we write
$\sigma(i)\;\mbox{for}\;\sigma_{\alpha}(i)$. However, sometiems we expilicitly write 
$\sigma_{\alpha}(i)$ to indicate to the reader the dependence of $\sigma(i)\;\mbox{on}\;
\alpha$.   
%\end{eqnarray*}
%Given a simple random walk $S_{i}$, we define
%the $N^{\alpha}$ loop erased walk by erasing $N^{\alpha}$ loops from the
%paths of $S_{i}$. We denote this process by $\hat{S}_{i}^{(N)}$. We omit
%$\alpha$ to avoid cumbersome notations.
%The details of the loop erasing procedures used in the construction of
%the $N^{\alpha}$ loop erased walk can be found in [L5].\\
\noindent  Our goal is to find $\lim_{N\rightarrow\infty}\frac{{\hat S_{N}}^{(N)}}{N^{\gamma}},\;
\mbox{for some}\;\gamma\;\mbox{in}\; Z^{3}.$
%\begin{equation}
%\lim_{N\rightarrow\infty}\frac{{\hat S_{N}}^{(N)}}{N^{\gamma}},\;
%\mbox{for some}\;\gamma\; in Z^{3}
%\end{equation}
We say that $n$ belongs to an $N^{\alpha}$ loop if $\exists
i\;\mbox{and}\;j\; \mbox{such that}\; i\leq n\leq
j\;\mbox{with}\;S_{i}=S_{j}\;\mbox{and}\;|i-j|\leq N^{\alpha}.$
For each $n$ we say $n$ is $N^{\alpha}$ loop free if $n$ does not belong
to an $N^{\alpha}\;\mbox{loop}.$ Suppose $n$ is $N^{\alpha}$ loop free. Then
loop erasing before $n$ and after $n$ are independent. If $n$ is
$N^{\alpha}$ loop free, then $n$ is not erased. However the converse is
not in general true. In order to analyze the behavior of $\hat{S}_{N}^{(N)}$ for
large $N$ we need to investigate how many steps of the SRW remain after
the first $N^{\alpha}$ loops have been erased.  Note that we  may still
have some small loops remaining after the first $N^{\alpha}$ loops  have
been erased. However, the algorithm to generate ${\hat S}_{n}$ only 
requires finite memory depending on $n$.  
%Using similar approach as in [L5], we introduce the following notions. The 
%notations
%we use are similar to the notations used in [L5].\\
%\noindent Let
%\begin{eqnarray*}
%$\sigma_{\alpha}(0)=\mbox{sup}\{j: S(j)= 0, |j|\leq\;N^{\alpha}\}$, and
%for $i>0\;\\
%\sigma_{\alpha}(i)=\mbox{sup}\{j>\sigma_{i-1}:
%S(j)=S(\sigma_{i-1}+1), |j-\sigma_{i-1}-1|\leq\;N^{\alpha}\}$.
%\end{eqnarray*}
Let $\rho_{\alpha}(j)=i\;\mbox{if}\;\sigma_{i}\leq j<\sigma_{i+1}.$
Then,$\rho_{\alpha}(\sigma(i))=i,\;\sigma(\rho_{\alpha}(j))\leq j$.
%\;\mbox{and}\;{\hat{S}_{i}^{(N)}}\;=\;S(\sigma(i)).
%\end{eqnarray*}
%\noindent Let $Y_{n}$ be the indicator function of the event "the
%$n^{th}$ point of S is not erased", i.e.
\noindent Let $Y_{n}=1\;\mbox{if}\;\sigma(i)=n\;\mbox{for some}\; 
i\geq 0,\;\mbox{and}\;Y_{n}=0\;\mbox{otherwise}$. Then 
$\rho(n)=\sum_{j=0}^{n}Y_{j}$ is the number of points remaining of the first $n$
points after the first $N^{\alpha}$ loops are erased. Let 
$a_{n,\alpha}=E(Y_{n,\alpha})$ be the probability that the $n^{th}$ point
is not erased. For the asymptotic behavior of $\rho_{\alpha}(N)$, we
have, \\
%\begin{thm} For $d=3,\;\mbox{and}\;0< \alpha <
%h(\zeta),\;\mbox{where}\;h(\zeta)=\frac{1}{1+2\zeta}$\\
%$$\frac{\rho_{\alpha}(N)}{Na_{N,\alpha}}\rightarrow\;1,$$ in probability
%as
%$N\rightarrow \infty.$
%\end{thm}
\vspace{2mm}

\noindent {\bf Theorem 2.1} For $ 0 < \alpha <\frac{1}{1+2\zeta},\;
\frac{\rho_{\alpha}(N)}{Na_{N,\alpha}}\rightarrow 1$ in probability as
$N\rightarrow \infty.$\\

\noindent G. Lawler proved analogous results in higher dimensions for
$\alpha=\infty$ (see [L1]).
%\begin{thm}
%a) For $d \geq 5,$ there exists an a $>0$ such that with probability
%one, $$\frac{\rho(N)}{N}\rightarrow\;\mbox{a},$$ as $N\rightarrow
%\infty.$\\\noindent
%b) For $d=4,$ $$\frac{\rho(N)}{Na_{N}}\rightarrow\;1,$$ in probability
%as $N\rightarrow \infty.$\\\\\noindent
%c) If $d \geq 5,$ and ${\tilde W}_{n}(t)=\frac{d \sqrt{a}
%\tilde{S}([nt])}{\sqrt n}$, then $\tilde{W}_{n}(t)\Rightarrow B(t)$,
%where $B(t)$ is a standard Brownian motion.\\\\\noindent
%d) If $d = 4,$ and ${\tilde W}_{n}(t)=\frac{d \sqrt{a_{n}}
%\tilde{S}([nt])}{\sqrt n}$, then $\tilde{W}_{n}(t)\Rightarrow B(t)$,
%where $B(t)$ is a standard Brownian motion.
%\end{thm}
%Here $\Rightarrow$ denotes weak convergence in the metric space
%$C([0,1])$ with the sup norm. 
%Theorem 2.1 says that if the length of the first $N^{\alpha}$ loops that
%are erased is less than $N^{h(\zeta)}$, then the number of points
%remaining is about $Na_{N,\alpha}$. Using this theorem we can then prove
%that, under the assumption $0 < \alpha < h(\zeta)
%\;\mbox{and}\; a_{N,\alpha}\sim \frac{\mbox{const.}}{N^{q}},\;\mbox{for
%some}\; q > 0 $, the limiting distribution of the $N^{\alpha}$ LEW is
%Gaussian.\\\\
\noindent Our next result is about the limiting distribution of the
$N^{\alpha}$ LEW. Let
$F_{N}$ be  defined by
$F_{N}\;=\;[{\sigma_{\alpha}(N)a_{\sigma_{\alpha}(N)}}].$ Here by
$\lceil\cdot \rceil$ we mean the greatest integer function. 
%Thus $Y_{N}$ is a positive integer valued random variable. 
Then we have,\\

\noindent {\bf Theorem 2.2}
(a) $\frac{S_{F_{N}}}{\sqrt N} \rightarrow \Phi, \;\mbox{where}\;
\Phi\;\mbox{is a normal random variable.}$\\
(b) Suppose $a_{N,\alpha}\sim \mbox{const}\cdot N^{-q},\; \mbox{for
some}\;q > 0.\;\mbox{Let}\; \tau_{N}=N^{-q/(1-q)}$. Then $
\frac{S_{\sigma_{\alpha}(N)}\sqrt \tau_{N}}{\sqrt N}\rightarrow \Phi.$\\
%\begin{thm}
%(a) $\frac{S_{Y_{N}}}{\sqrt N} \rightarrow \Phi, \;\mbox{where}\;
%\Phi\;\mbox{is a normal random variable.}$\\
%(b) Suppose $a_{N,\alpha}\sim \mbox{const}\cdot N^{-q},\; \mbox{for
%some}\;q > 0.\;\mbox{Let}\; \tau_{N}=N^{-q/(1-q)}$. Then $$
%\frac{S_{\sigma_{\alpha}(N)}\sqrt \tau_{N}}{\sqrt N}\rightarrow \Phi.$$
%\end{thm}

\noindent Clearly, $q$ satisfies $0 < q \leq \alpha \zeta.$ However, we
were unable to prove the existence of $q$.
%Theorem 2.2 shows that,in $Z^{3}$, the behavior of the $N^{\alpha}$ LEW
%for $\alpha < h(\zeta)$ is the same as that of the SRW. 
For a sufficiently large $\alpha$ we have,\\
%the following result.
%can
%conclude that the limiting  distribution of the
%$N^{\alpha}$ LEW is the same as the limiting  distribution
%of the LEW. In fact we have,\\
\vspace{1mm} 
\noindent {\bf Theorem 2.3}
Let $c_{N}=(E(|S_{\sigma_{\alpha}(N)}|^{2}))^{1/2},\;\mbox{and}\;
d_{N}=(E(|S_{\sigma_{(N)}}|^{2}))^{1/2}.$ Suppose that
$\frac{S_{\sigma_{\alpha}(N)}}{c_{N}}\;\mbox{or}\;\frac{S_{\sigma{(N)}}}{d_{
N}}$ converge in distribution. If $\alpha > 2$, then
$\lim_{N\rightarrow \infty}
\frac{S_{\sigma_{\alpha}(N)}}{c_{N}}=\lim_{N\rightarrow
\infty}\frac{S_{\sigma(N)}}{d_{N}},$ in distribution. Here
$\sigma(N)=\sigma_{\alpha}(N)\;\mbox{with}\;\alpha=\infty.$

%\noindent The limiting distribution of the loop erased walk in three
%dimensions is still
%unknown. It is believed, however,  that it is non-Gaussian (see [L5]).
%In this
%direction we have the following result.\\
%\begin{thm} Let$\;$
%$c_{N}=(E(|S_{\sigma_{\alpha}(N)}|^{2}))^{1/2},\;\;\mbox{and}\;\;
%d_{N}=(E(|S_{\sigma_{(N)}}|^{2}))^{1/2}.$\\\noindent  Suppose that
%$\frac{S_{\sigma_{\alpha}(N)}}{c_{N}}\;\mbox{or}\;\frac{S_{\sigma{(N)}}}{d_{
%N}}$ converge in distribution.\\\noindent  If $\alpha > 2$, then
%$$\lim_{N\rightarrow \infty}
%\frac{S_{\sigma_{\alpha}(N)}}{c_{N}}=\lim_{N\rightarrow
%\infty}\frac{S_{\sigma(N)}}{d_{N}},$$ in distribution. Here
%$\sigma(N)=\sigma_{\alpha}(N)\;\mbox{with}\;\alpha=\infty.$
%\end{thm}
%Here $\sigma(N)=\sigma_{\alpha}(N)\;\mbox{with}\;\alpha=\infty.$ 
%Hence the limiting distribution of the $N^{\alpha}$ loop erased walk is
%equal to the limiting distribution of the loop erased walk in three
%dimensions for sufficiently big $\alpha$. 
%This theorem gives a theoretical
%justification for doing computer simulation of the loop erased walk in
%three dimensions. Since the limiting distribution of the $N^{\alpha}$
%loop erased walk for $\alpha > 2$ is equal to the limiting distribution
%of the loop erased walk it would be possible to simulate the loop erased
%walk by erasing loops whose length is at most $N^2$. That way one can
%determine from a numerical point of view what the limiting distribution
%of the loop erased walk is in three dimensions.
\vspace{3mm}

\noindent {\bf 3. PROOFS}\\

\noindent $\mbox{For}\;0\leq j< k<\infty$, we denote by $Z(j,k)$
the indicator function of the event "there is no $N^{\alpha}$ loop free
point between $j\;\mbox{and}\;k$ including $j\;\mbox{and}\;k$".\\

\noindent {\bf Lemma 3.1} There exist constants $c_{1}, c_{2}$ such that if $\beta > \alpha$, 
then $ E(Z(k-N^{\beta},k))\leq c_{1}e^{-c_{2}N^{\beta-\alpha}}.$\\

\noindent {\bf Proof:} From Theorem 1.1 of [L2] it follows that there is a 
$c_{3}$ such that in the interval $[k-4N^{\alpha}]$ the probability of an 
$N^{\alpha}$ loop free point is at least $c_{3}$. Consider now an interval I of 
length $N^{\beta}$ divided into $\frac{1}{4}N^{\beta-\alpha}$ small
intervals of length $4N^{\alpha}$. Then the probability of no $N^{\alpha}$
loop free point in I is bounded by
$(1-c_{3})^{\frac{1}{4}N^{\beta-\alpha}}$ which can be written  in the
form $c_{1}e^{-c_{2}N^{\beta-\alpha}}$.\\

\noindent Suppose that for some $k,\;\mbox{where}\;0\leq k\leq N,\;
N^{\alpha}$ loops are erased only on $S[k,\infty)$, so that $S_{k}$ is 
considered to be the origin. Let $Y_{N,k}$ be the probability that $S_{N}$
is not erased in this procedure. Clearly $E(Y_{N,k})=a_{N-k}.$
Now suppose $0\;\leq\;k\;\leq N-N^{\beta},\; \mbox{for some}\; \beta
< 1\; \mbox{and}\; Z(N-N^{\beta}, N)=0$. 
%i.e. that there exists an $N^{\alpha}\;\mbox{loop free point
%between}\;N-N^{\beta}\;\mbox{and}\; N$. 
Then it can be shown that
$Y_{N,k}=Y_{N}$, and hence by Lemma 3.1,
%\begin{eqnarray*}
 $|a_{N}-a_{N,k}|\;\leq \;P\{Y_{N}\neq Y_{N,k}\} \leq 
E(Z(N-N^{\beta}, N))\leq c_{1}e^{-c_{2}N^{\beta-\alpha}}.$
%\end{eqnarray*}
Thus, for $N^{\beta}\;\leq\;k\;\leq N$,
\begin{equation}
\noindent |a_{k}-a_{N}|\;\leq\;c_{1}e^{-c_{2}N^{\beta-\alpha}}\leq 
c_{1}a_{N}N^{\alpha\zeta}e^{-c_{2}N^{\beta-\alpha}}
\end{equation}
\noindent {\bf Proof of Thm 2.1} For each $N$, choose
$0\leq j_{0}<j_{1}<j_{2}<...<j_{m}=N,$ such that $j_{i}-j_{i-1}\sim
N^{1-\alpha\zeta-\delta},$ uniformly in i.
Then
$m\;\sim\;N^{\alpha\;\zeta+\delta}.\;$ Erase loops on each interval
$[j_{i},j_{i+1}]$ separately. Let $\tilde{Y}_{k}$ be the indicator
function of the event "$S_{k}$ is not erased in this finite
loop-erasing". Let $K_{0}=[0,0],$ and $\epsilon_{1}>\delta$. Then, for
$i=1,...,m,\;$ define the intervals $K_{i}\;\mbox{and}\;K_{i}^{'}$ by
$K_{i}=[j_{i}-N^{1-2\alpha\zeta-\epsilon_{1}},j_{i}],\;\;
K_{i}^{'}=[j_{i},j_{i}+N^{1-2\alpha\zeta-\epsilon_{1}}].$ Let
$R_{i},\;  i=1,...,m,$ be the indicator function of the event  
%of the complement of the
%event,"there exist $N^{\alpha}$ loop free points in both $K_{i}^{'}$ and
%$K_{i+1}$", i.e.
$\{\exists{\mbox{ no}\;N^{\alpha}\;\mbox{loop free point
in}\;K_{i}^{'}\; \mbox{or in}\; K_{i+1}\}}.$
Note that $R_{i}=0 \;\mbox{if and only if} \; \exists
N^{\alpha}\;\mbox{loop free point in}\;K_{i}^{'}\;\mbox{and in}\;K_{i+1}.$
Thus if $j_{i}+N^{1-2\alpha\zeta-\epsilon_{1}}\leq k \leq
j_{i+1}-N^{1-2\alpha\zeta-\epsilon_{1}}\;\mbox{and}\;
R_{i}\;=\;0,$ then $Y_{k}=\tilde{Y_{k}}.$ Therefore for a sufficiently large 
$N$,
%\bigskip
\begin{equation}
|\sum_{k}\;Y_{k}-\tilde{Y_{k}}|\;\leq
2N^{1-\alpha\zeta-\epsilon_{1}+\delta}+ 2N^{1-\alpha\zeta-\delta}
\sum_{i}\;R_{i}.
\end{equation}

%2N^{\alpha\zeta+\delta}\frac{N}{N^{2\alpha\zeta+\epsilon_{1}}}+\frac{2N}
%{N^{\alpha\zeta+\delta}}\sum_{i}\;R_{i}\nonumber \\ &=&

\noindent Let $\lambda=1-2\alpha\zeta-\epsilon_{1}-\alpha$. Then, 
\begin{eqnarray}
P\{\sum_{i}\;R_{i}\;\geq\;\frac{1}{4} N^{\gamma}\}&\leq&
4 c_{1}e^{-c_{2}{N^{\lambda}}}N^{\alpha\zeta+\delta-\gamma}
\end{eqnarray}
Since $\epsilon_{1}$ is arbitrary, for $\alpha < \frac{1}{1+2\zeta}$,
$\lambda > 0$ and the right side of $(3)$ goes to $0$ as $N \rightarrow
\infty$. Let now $\epsilon_{2}<<
\mbox{min}\{\epsilon_{1}-\delta;\frac{\delta}{2}\}$. Then using $(2)$ we
get 
%we conclude that for any

%$P\{\sum_{i}\;R_{i}\;\geq\;\frac{1}{4} N^{\gamma}\} \rightarrow 0$ as $N 
%\rightarrow \infty$.\\\noindent  
 
%$$P\{\sum_{k}\;Y_{k}-\tilde{Y}_{k}\geq\;\frac{N}{N^{\alpha\zeta+\epsilon_{2}}}\
%}.
\begin{eqnarray}
P\{\sum_{k}\;Y_{k}-\tilde{Y}_{k}\geq\;N^{1-\alpha\zeta-\epsilon_{2}}\}
\leq
%P\{\frac{2N}{N^{\alpha\zeta+\epsilon_{1}-\delta}}+\frac{2N}{N^{\alpha\zeta+
%\delta}}\sum\;R_{i}\geq
%\frac{N}{N^{\alpha\zeta+\epsilon_{2}}}\}\nonumber  
%\\&=&
%P\{\frac{2}{N^{\alpha\zeta+\epsilon_{1}-\delta}}+\frac{2}{N^{\alpha\zeta+\delta
%}}\sum\;R_{i}\geq \frac{1}{N^{\alpha\zeta+\epsilon_{2}}}\}\nonumber
%\\&\leq& P\{\frac{2}{N^{\alpha\zeta+\delta}}\sum 
%R_{i}\geq\frac{1}{2N^{\alpha\zeta+\epsilon_{2}}}\} 
P\{\sum_{i} R_{i} \geq \frac{1}{4}N^{\delta-\epsilon_{2}}\}.
\end{eqnarray}
Put $\delta-\epsilon_{2}=\gamma$.
Then $(4)$ goes to $0$ by $(3)$. From
$(4)$ it  follows that
%$$\frac{1}{Na_{N}}\sum (Y_{k}-\tilde{Y_{k}}).$$
%\begin{eqnarray}
%\noindent P\{\frac{1}{Na_{N}}\sum Y_{k}-\tilde{Y_{k}}\geq \tilde{\epsilon}\}&=&
%P\{\sum Y_{k}-\tilde{Y_{k}}\geq\tilde{\epsilon} Na_{N}\}\nonumber \\
%&\leq&
%P\{\sum
%Y_{k}-\tilde{Y_{k}}\geq \tilde{\epsilon} \frac{N}{N^{\alpha\zeta}}\}.
%\end{eqnarray}
%\bigskip
%Set
%\begin{center}
%$\frac{\tilde{\epsilon}N}{N^{\alpha\zeta}}=\frac{N}{N^{\alpha\zeta+\epsilon_{2}
%}}\;\mbox{and choose}\; \tilde{\epsilon}\;\mbox{to 
%be}\;\frac{1}{N^{\epsilon_{2}}}.$
%\end{center}
%Using $(8)\;\mbox{and}\;(9)$ we have, thus , proved
$\frac{1}{Na_{N}}\;\sum_{k}\;Y_{k}-\tilde{Y_{k}}\rightarrow\;0\;\mbox{in probability}.$
We can write $\sum \tilde{Y_{k}}=1+\sum X_{i},$ where $X_{i}$ are the independent random 
variables, $X_{i}=\sum_{k=j_{i-1}}^{{j_i}-1}\;\tilde{Y_{k}}.$ 
Then, using (3) and Chebyshev's Inequality, we can show, 
$\frac{1}{E(\sum_{k} \tilde Y_{k})}\sum_{k} \tilde Y_{k} \rightarrow
1\;\mbox{in probability}$. From $(3)$ and Lemma 3.1 follows that
%\end{equation}. From $(3)$ and Lemma 3.1 follows that
$E(\sum_{k=0}^{N}
\tilde{Y}_{k}) \sim Na_{N}$, completing the proof of the theorem. 
\vspace{3mm}
\noindent {\bf Propositon 3.1} Let $\sigma(N)=\sigma_{\alpha}(N)$ be defined as in section
$2$. Then\\\noindent  (a) $\frac{\sigma(N)a_{\sigma(N)}}{N}
\rightarrow 1\;\mbox{in probability as}\;N \rightarrow \infty.$\\
(b) Assume $a_{N} \sim \frac{1}{N^{q}}\; \mbox{for some}\; q >
0\;\mbox{and let}\;\tau_{M}\sim M^{-q/(1-q)}.$\\
Then $\frac{\sigma(M)\tau_{M}}{M}\rightarrow 1\;\mbox{in probability as}\; M
\rightarrow \infty.$\\

\noindent {\bf Proof of (a):} Let $s>0$ be a constant. It suffices to prove that
$\frac{\sigma(M_{t})a_{\sigma(M_{t})}}{M_{t}}$ converges to $ 1\;\mbox{a.s.for 
any sequence}\; M_{t}\geq t^{s}.$ 
%$||\frac{\rho_{\alpha}(N)}{Na_{N}}-1||\leq\;\frac{C}{N^{\delta}},\;\mbox{f
%or some}\; \delta > 0,\;\mbox{we
%have}\;\frac{\rho_{\alpha}(N_{t})}{N_{t}a_{N_{t}}}-1
%\rightarrow
%0\; a.s.\;\mbox{for}\; N_{t}=t^{s}\;\mbox{with}\; s\delta > 1.$
%Therefore,
%\begin{equation}
By Theorem 2.1 $\exists\; \Omega^{'}\subset \Omega\\\mbox{such that}\;
P(\Omega^{'})=1$ and $\frac{\rho_{\alpha}(N_{t})}{N_{t}a_{N_{t}}}
-1 \rightarrow 0,\; \forall \omega \in \Omega^{'}$. Let $N_{t}^{'}$ be a
%\end{equation} 
sequence such that
$N_{t}^{'}\geq t^{s}$.
%$\rho_{\alpha}(N_{t}^{'})(\omega)(N_{t}^{'}a_{N_{t}^{'}})^{-1}
%\rightarrow 1 
%\forall \omega \in \Omega^{'}$. 
Then for a fixed $t$ there exists a sequence $\xi_{t}$ such that
$\;\xi_{t}^{s}\leq N_{t}^{'}< (\xi_{t}+1)^{s}$. Note that $t \leq
\xi_{t}.$ For $\omega \in \Omega^{'},$ 
\begin{equation}
\rho_{\alpha}(\xi_{t}^{s})((\xi_{t}+1)^{s}a_{N_{t}^{'}})^{-1}\leq
\rho_{\alpha}(N_{t}^{'})(N_{t}^{'}a_{N_{t}^{'}})^{-1} \leq
\rho_{\alpha}(\xi_{t}+1)^{s}(\xi_{t}^{s}a_{N_{t}^{'}})^{-1}.
\end{equation}
By Theorem 2..1 and (2) the upper and lower bounds of this inequality
converge to 1 in probability. Substituting $\sigma(M_{t})$ for
$N^{'}_{t}$ gives
$M_{t}(\sigma(M_{t})a_{\sigma(M_{t})})^{-1}\rightarrow 1$\\ 

\noindent {\bf{Proof of (b)}} From (a) we have 
%$a_{\sigma_{M_{t}}(\omega)}(\tau_{M_{t}})^{-1}\rightarrow 1\; \forall
%\omega \in \Omega^{'}$. From (a), we have,
$\sigma(M_{t})a_{\sigma(M_{t})}(M_{t})^{-1}\rightarrow 1.$
By assumption,
$\frac{\sigma(M_{t})\sigma(M_{t})^{-q}}{M_{t}}\rightarrow 1.$
Therefore,
%$\frac{\sigma(M_{t})^{1-q}(\omega)}{M_{t}}\rightarrow 1,$
%or,
$\frac{\sigma(M_{t})(\omega)}{M_{t}^{1/(1-q)}}\rightarrow 1, \mbox{as
t}\;\rightarrow \infty.$
Since this holds for all $M_{t}\geq t^{s},\;
\sigma(N)(N^{1/(1-q)})^{-1}\rightarrow 1\;\mbox{in probability}$. By Proposition 3.1a,
$\frac{{N^{1/(1-q)}}}{\sigma(N)}\cdot \frac{\sigma(N)a_{\sigma(N)}}{N}\rightarrow 1\;
\mbox{in probability}.$ This and Proposition
3.1a  imply $a_{\sigma(N)}(\tau_{N})^{-1}\rightarrow 1\;\mbox{in probability}.$ 
Using Prop.3.1a again, we get, $\sigma(N)\tau_{N}(N)^{-1}\rightarrow 1\;\mbox{in probability}.$
Hence
$[\sigma_{\alpha}(N)a_{\sigma_{\alpha}(N)}](N)^{-1}\rightarrow 1\;\mbox{in probability}$.\\
% \noindent  Let $Y_{N}$ be defined by
%$Y_{N}\;=\;[{\sigma_{\alpha}(N)a_{\sigma_{\alpha}(N)}}]$.\\

\noindent {\bf Proof of Thm 2.2} The proof of Theorem 2.2 follows
from Theorem 3.1 and Proposition 3.1.\\\noindent  
%Statement b of
%Theorem 2.3 follows from  Theorem 3.1 and  Proposition 3.1(b).

\vspace{3mm}
\noindent {\bf Theorem 3.1} Let $X_{j}$ be i.i.d. random variables with
$E(X_{j})=0$ and $Var(X_{j})=1.$ Let
$\nu_{n}$ be positive integer valued random variables such that
$\frac{\nu_{n}}{n}\rightarrow c\;\mbox{in probability}$. Then
$\frac{S_{\nu_{n}}}{\sqrt{cn^{p}}}$ converges in distribution to a
standard normal random variable $\cal{N}.$
%\mbox{with}\;
%\cal{N}\;\mbox{is a normal random}\\\noindent
%\mbox{variable with mean zero and has the identity matrix variance}.$
\vspace{3mm}

\noindent {\bf Proof of Thm 2.3}$\;$For each $N,$ choose $0\leq
j_{1}<j_{2}...<j_{m},$ satisfying
$j_{i}-j_{i-1}\sim\;N^{\alpha},N^{s}-N^{t}\;\leq j_{i}\leq\;N^{s}$.
Let $ X=\sum_{i=1}^{m}1_{\{j_i\}}.$ Then,
$E(Z(N^{s}-N^{t}, N^{s}))\leq c_{1}e^{-c_{2}N^{t-\alpha}}.$  
%P\{X=0\}\nonumber \\&\leq& 
%P\{X\leq \frac{1}{2}E(X))\nonumber \\&=&P(|X-E(X)|>\frac{1}{2}E(X))\nonumber
%\\&\leq&\frac{4}{E(X^{2})}Var(X).
%$$Var(X)=\sum_{i=1}^{m}Var(1_{\{j_{i}\}})\;\leq\;\sum_{i=1}^{m}E(1_{\{j_{i}\}})
%\=\;E(X).$$
%Now,
%\begin{eqnarray*}
%E(Z(N^{s}-N^{t}, N^{s}))\;&\leq&\frac{4}{E(X^{2})}Var(X)\nonumber
%\\&\leq&
%\frac{4}{E(X)}.
%\end{eqnarray*}
%But,
%\begin{eqnarray*}
%E(X)\;&\geq&\;N^{t-\alpha}P\{N^{s}\; \mbox{is}\; N^{\alpha}\; \mbox{loop 
%free}\}\nonumber \\&\geq&\; N^{t-\alpha}N^{-\alpha \zeta}\nonumber \\&=& 
%N^{t-\alpha-\alpha \zeta}.
%\end{eqnarray*}
%We want $$t-\alpha-\alpha \zeta\; >\;0.$$
%Since  $$\alpha+\alpha \zeta\;<\;\alpha+\frac{\alpha}{2},$$ choosing
%$t\;>\;\frac{3\alpha}{2}$ guarantees the positivity of the term
%$t-\alpha-\alpha
%\zeta.$\\
Consider the interval $[0,N^{s}]$ divided into subintervals of
length $N^{t},t<s$. Then,
%\begin{eqnarray*}
$P\{\rho(N^{s})\;<\;N^{s-t}\}
%&\leq& N^{s-t}P\{[N^{s}-N^{t}, N^{s}]\;
%\mbox{has no loop free point}\}\nonumber \\ &\leq&
\leq c_{1}N^{s-t}e^{-c_{2}N^{t-\alpha}}$
%&\leq& N^{s-2t+3\alpha/2}.
%\end{eqnarray*}
and
%\begin{eqnarray}
$P\{N^{s}<\;\sigma(N^{s-t})\} \leq
P\{\rho(N^{s})\;<\;\rho(\sigma(N^{s-t}))\}
%\nonumber \\ &=& P\{\rho(N^{s})< N^{s-t}\}\nonumber \\
\leq c_{1}N^{s-t}e^{-c_{2}N^{t-\alpha}}.$
%\end{eqnarray}
Let $M=N^{s-t}$. Then,
%\begin{eqnarray}
$P\{\sigma(M)>M^{\frac{s}{s-t}}\}\leq (c_{1}M)e^{-c_{2}M^{\frac{t-\alpha}{s-t}}}.$
%\\&=& M^{\frac{s-2t+3\alpha/2}{s-t}}.
%\end{eqnarray}
We show the $L^{2}$ norm
of the difference of the $N^{\alpha}$ LEW and the LEW goes to $0$. Let
$e_{N}=\mbox{max}\{c_{N},d_{N}\}.$ Then, 
%$$||\frac{S_{\sigma_{\alpha}(N)}}{c_{N}}-\frac{S_{\sigma(N)}}{d_{N}}||_{2}\;
%\leq\;
%||\frac{S_{\sigma_{\alpha}(N)}}{d_{N}}-\frac{S_{\sigma
%(N)}}{d_{N}}||_{2}+
%||\frac
%{S_{\sigma_{\alpha}(N)}}{c_{N}}-\frac{S_{\sigma_{\alpha}(N)}}{d_{N}}||_{2}.$$
%Consider the term:
%$$||\frac{S_{\sigma_{\alpha}(N)}}{c_{N}}-\frac{S_{\sigma_{\alpha}(N)}}
%{d_{N}}||_{2}.$$
%We then have,
%\begin{eqnarray}
%||\frac{S_{\sigma_{\alpha}(N)}}{c_{N}}-\frac{S_{\sigma_{\alpha}(N)}}{d_{N}}|
%|_{2}&=&||\frac{S_{\sigma_{\alpha}(N)}}{c_{N}}||_{2}\;
%|\frac{d_{N}-c_{N}}{d_{N}}|\nonumber
%\\&=&|\;\frac{||S_{\sigma(N)}||_{2}-||S_{\sigma_{\alpha}(N)}||_{2}}{d_{N}}|\\
%&\leq&\frac{||S_{\sigma(N)}-S_{\sigma_{\alpha}(N)}||_{2}}{d_{N}}.
%\end{eqnarray}
%$(30)$ follows since $||\frac{S_{\sigma_{\alpha}}(N)}{C_{N}}||\;=\;1\;
%\mbox{by definition}.$\\
%Similarly,
%$$||\frac{S_{\sigma_{\alpha}(N)}}{c_{N}}-\frac{S_{\sigma 
%(N)}}{d_{N}}||_{2}\;\leq\;||
%\frac{S_{\sigma_{\alpha (N)}}}{c_{N}}-\frac{S_{\sigma
%(N)}}{c_{N}}||_{2}+
%||\frac{S_{\sigma(N)}}{c_{N}}-\frac{S_{\sigma(N)}}{d_{N}}||_{2}.$$
%Consider the term
%\begin{eqnarray}
%||\frac{S_{\sigma(N)}}{c_{N}}-\frac{S_{\sigma(N)}}{d_{N}}||_{2}&=&||
%\frac{S_{\sigma(N)}}{d_{N}}||_{2}\; |\frac{c_{N}-d_{N}}{c_{N}}|\nonumber
%\\&=&|\;\frac{||S_{\sigma_{\alpha}}(N)||_{2}-||S_{\sigma}(N)||_{2}}{c_{N}}|\
%\&\leq&\frac{||S_{\sigma_{\alpha}(N)}-S_{\sigma(N)}||_{2}}{c_{N}}.
%\end{eqnarray}
%$(32)$ follows since $||\frac{S_{\sigma(N)}}{d_{N}}||\;=\;1\;
%\mbox{by definition}.$\\
%and
%\begin{equation}
$||\frac{S_{\sigma_{\alpha}(N)}}{c_{N}}-\frac{S_{\sigma(N)}}{d_{N}}||_{2}\leq
\frac{2 \cdot ||S_{\sigma(N)}-S_{\sigma_{\alpha}(N)}||_{2}}{e_{N}}.$\\
%\end{equation}
\noindent Let ${\cal S} = \{\omega \in \Omega: \exists \mbox{a loop
between}\;i\;\mbox{and}\;j\;0\leq i \leq \sigma_{\alpha(N)}, 
|i-j|>N^{\alpha}\}$, ${\cal T} =
\{\omega \in \Omega: \exists \mbox{a loop between}\;i\;\mbox{and}\;j\;
0\leq i \leq \sigma_{\alpha(N)},\;\sigma (N) > N^{\frac{s}{s-t}},\;
|i-j|>N^{\alpha}\}$ and 
$G=\{\omega \in \Omega:|\frac{S_{\sigma(N)}-S_{\sigma_{\alpha}(N)}}{e_{N}}|>M\}.$ Let 
${\cal I}$ be the indicator function defined on ${\cal S}$.
%Then ${\cal I}=1_{\cal T}+ 1_{{\cal T}^{c}}}$
%of the event "there exists a loop between
%$i\;
%\mbox{and}\; j$ with loop length greater than
%$N^{\alpha}$".
%i.e.\\
%${\cal I}=1_{\{\exists \mbox{a loop between}\;i\;\mbox{and}\; j,\;
%0\leq i \leq \sigma_{\alpha(N\}}, |i-j|>N^{\alpha}\}}$. Let $G=\{\omega
%\in 
%\Omega:|\frac{{S_{\sigma(N)}-S_{\sigma_{\alpha}(N)}}}{e_{N}}|>M\}.$\\
\noindent
Then,
%\begin{eqnarray}
$\frac{||S_{\sigma(N)}-S_{\sigma_{\alpha}(N)}||_{2}^{2}}{e_{N}^{2}}=
E(\frac{|S_{\sigma(N)}-S_{\sigma_{\alpha}(N)}|^{2}}{e_{N}^{2}}\cdot {\cal 
I})=\\
\int_{G}\frac{|S_{\sigma(N)}-S_{\sigma_{\alpha}(N)}|^{2}}{e_{N}^{2}}\cdot{\cal 
I} dP +
\int_{G^{c}}\frac{|S_{\sigma(N)}-S_{\sigma_{\alpha}(N)}|^{2}}{e_{N}^{2}}\cdot
{\cal I}dP.$
%\end{eqnarray}
Then $\forall \epsilon > 0 \;\exists\; M_{0}\; \mbox{such 
that}\;\\ \forall M \geq M_{0},
%\nonumber \\&\leq&
%Since
%\begin{eqnarray}
%(\int_{\Omega}\frac
%{|S_{\sigma(N)}-S_{\sigma_{\alpha}(N)}|^{2}}{e_{N}^{2}})^{1/2}&\leq&
%\frac{1}{e_{N}}(||S_{\sigma(N)}||_{2}+||S_{\sigma_{\alpha}(N)}||_{2})\nonumber\
%\
%&\leq&\frac{||S_{\sigma(N)}||_{2}}{d_{N}}+\frac{||S_{\sigma_{\alpha}(N)}||_{2}}
%{c_{N}}
%\nonumber \\ &=& 2\; < \infty,
%\end{eqnarray}
%by Dominated Convergence Theorem,\\
%\noindent $\forall \epsilon > 0 \;\exists M_{0}\;
%\mbox{such that}\; \forall M \geq M_{0},$\\
%\begin{eqnarray}
\int_{G}\frac{|S_{\sigma(N)}-S_{\sigma_{\alpha}(N)}|^{2}}{e_{N}^{2}}\cdot{\cal
I} dP 
%\int_{G}\frac{|S_{\sigma(N)}-S_{\sigma_{\alpha}(N)}|^{2}}{e_{N}^{2}}dP
< \;\epsilon, \; \mbox{if}\; M \geq M_{0}.$
%\end{eqnarray}
Consider now the second summand with $M=M_{0}$.
%\pagebreak
%\begin{eqnarray}
$\int_{G^{c}}\frac{|S_{\sigma(N)}-S_{\sigma_{\alpha}(N)}|^{2}}{e_{N}^{2}}\cdot
{\cal I}dP \leq M_{0}^{2}\int_{\Omega}{\cal I} dp
= M_{0}^{2}\cdot {E1_{{{\cal T}^{c}}}}
+ E1_{{\cal T}}
\leq M_{0}^{2} \cdot\{  E1_{{{\cal T}^{c}}}\}+ 
c_{1}N^{s-t}e^{-c_{2}N^{t-\alpha}}\} \sim M_{0}^{2}\cdot
\sum_{i=0}^{N^{\frac{s}{s-t}}}\sum_{j=N^{\alpha}}^{\infty}
\frac{1}{|j-i|^{3/2}}+N^{\frac{s-2t+3\alpha/2}{s-t}} \leq
M_{0}^{2}(\frac{N^{\frac{s}{s-t}}}{N^{\alpha/2}}+c_{1}N^{s-t}e^{-c_{2}N^{
t-\alpha}}).$
%\end{eqnarray}
%\noindent $(38)$ follows from $(29).$\\\noindent
\noindent For $\alpha > 2\; \mbox{there exist}\;s\;\mbox{and}\; t$ such that
the last term goes to $0$.\\

\noindent {\bf Acknowledgement} The authors are thankful to Professor Greg Lawler for
constructive comments on the previous version of this paper and for
suggesting a stronger version of Lemma 3.1.\\
\vspace{200mm}

\noindent {\bf REFERENCE}\\

\noindent [L1] Lawler, G. 1991. {\em Intersections of random walks.} (Birkh\"{a}user Boston).\\
\noindent [L2] Lawler, G., 1996. {\em Cut Times For Simple Random Walk,} EJP. {\bf Vol 1}: Paper 13.\\
\noindent [L3] Lawler, G., {\em Strict Concavity Of The Intersection
Exponent For Brownian Motion in 2 And 3 Dimensions}, Math Physics
Electronic Journal, 5 (1998).

\end{document}